\documentclass[12pt, a4paper, leqno]{amsart}
\usepackage[utf8]{inputenc}
\usepackage{amsmath}
\usepackage{amsfonts}
\usepackage{amssymb}
\usepackage{times}
\usepackage[T1]{fontenc} %skandit 
\usepackage{url} 
\usepackage{color,esint}
\usepackage[dvipsnames]{xcolor}
\usepackage[colorlinks, allcolors=RedViolet,pdfstartview=,pdfpagemode=UseNone]{hyperref} 

\usepackage{graphicx} 
\usepackage{enumitem}
\usepackage{tabularx}
    \usepackage{mathtools}
\usepackage{tikz}
\usepackage{pgfplots}
\pgfplotsset{compat=1.15}
\usepackage{mathrsfs}
\usetikzlibrary{arrows}
\usepackage{centernot} %Tämän avulla \imlies-komennon päälle tuleva kieltoviiva menee oikeaan paikkaan.

\let\oldmarginpar\marginpar
\renewcommand\marginpar[1]{\-\oldmarginpar[\raggedleft\footnotesize #1]%
{\raggedright\footnotesize #1}}

\theoremstyle{plain}
\newtheorem{theorem}[equation]{Theorem}%[section]
\newtheorem{lemma}[equation]{Lemma}
\newtheorem{proposition}[equation]{Proposition}
\newtheorem{corollary}[equation]{Corollary}

\theoremstyle{definition}
\newtheorem{definition}[equation]{Definition}

\theoremstyle{remark}
\newtheorem{remark}[equation]{Remark}

\numberwithin{equation}{section}

\newcommand{\R}{\mathbb{R}}

\newcommand{\N}{\mathbb{N}}
\newcommand{\Rn}{{\mathbb{R}^n}}

\renewcommand{\phi}{\varphi}
\renewcommand{\epsilon}{\varepsilon}
\def\le{\leqslant}

\def\ge{\geqslant}

\def\phi{\varphi}
\def\rho{\varrho}
\def\vartheta{\theta}

\newcommand{\Phiw}{\Phi_{\text{\rm w}}}

\def\supp{\operatorname{supp}}
\def\esssup{\operatornamewithlimits{ess\,sup}}
\def\osc{\operatornamewithlimits{osc}}

\def\dist{\qopname\relax o{dist}}
\def\osc{\operatornamewithlimits{osc}}
\def\essinf{\operatornamewithlimits{ess\,inf}}

\def\supp{\operatornamewithlimits{supp}}

\def\loc{{\rm loc}}

\newcommand{\ainc}[1]{\hyperref[ainc]{{\normalfont(aInc){\ensuremath{_{#1}}}}}}
\newcommand{\adec}[1]{\hyperref[adec]{{\normalfont(aDec){\ensuremath{_{#1}}}}}}
\newcommand{\inc}[1]{\hyperref[inc]{{\normalfont(Inc){\ensuremath{_{#1}}}}}}
\newcommand{\dec}[1]{\hyperref[dec]{{\normalfont(Dec){\ensuremath{_{#1}}}}}}
\newcommand{\azero}{\hyperref[azero]{{\normalfont(A0)}}}
\newcommand{\aone}{\hyperref[aonen]{{\normalfont(A1)}}}

\newcommand{\aonen}[1]{\hyperref[aonen]{{\normalfont(A1-\ensuremath{#1})}}}
\newcommand{\VAn}[1]{\hyperref[VA1n]{{\normalfont(VA1-\ensuremath{#1})}}}
\newcommand{\wVAn}[1]{\hyperref[wVA1n]{{\normalfont(wVA1-\ensuremath{#1})}}}
\newcommand{\condMn}[1]{\hyperref[condM]{{\normalfont(M-\ensuremath{#1})}}}

\date{\today}

\newcommand{\comment}[1]{\vskip.3cm
\fbox{%
\parbox{0.93\linewidth}{\footnotesize #1}}
\vskip.3cm}

\begin{document}

\title[Approximation results for non-standard growth functionals]{Optimal approximation results for non-standard growth functionals in the super-dimensional case}

\author{Peter Hästö}
\address{Peter Hästö, Department of Mathematics and Statistics,
FI-00014 University of Helsinki, Finland}
\email{\texttt{peter.hasto@helsinki.fi}}

\author{Jihoon Ok}
\address{Department of Mathematics, Institute for Mathematical and Data Science, Sogang University, Seoul 04107, Republic of Korea}
\email{jihoonok@sogang.ac.kr}


\begin{abstract}
In the super-dimensional case $p>n$ we show that $W^{1,p}$-functions can be uniformly approximated 
by Hölder continuous functions with improved exponent $\frac {p-n}{p-1}$. 
We use this to show density of smooth functions and absence of Lavrentiev's 
phenomenon for $(p, q)$-growth with optimal conditions. The special cases 
of double phase and Musielak--Orlicz growth are also included. 
Finally, we improve the conditions for maximal regularity of the double phase 
minimizers. In each case the optimality of our range of exponents is shown by known 
examples. 
\end{abstract}

\maketitle

%%%%%%%%%%%%%%%%%%%%%%%%%%%%%%%%%%%%%%%%%%%%%%%%%%%%%%%%%%%%%%%%%%%%%%%%%%%%%%%%%%%
%%%%%%%%%%%%%%%%%%%%%%%%%%%%%%%%%%%%%%%%%%%%%%%%%%%%%%%%%%%%%%%%%%%%%%%%%%%%%%%%%%%
%%%%%%%%%%%%%%%%%%%%%%%%%%%%%%%%%%%%%%%%%%%%%%%%%%%%%%%%%%%%%%%%%%%%%%%%%%%%%%%%%%%

\section{Introduction}

Approximation of Sobolev functions by smooth functions is a basic tool in the calculus of variations and the theory of partial differential equations.
For energies with nonstandard growth, however, a finite energy function need not admit smooth approximations whose energies converge. 
This failure is closely related to the Lavrentiev phenomenon: the infimum may depend on whether the admissible functions belong to a Sobolev space or a smooth class. 

A useful tool for these (and other) problems is truncating a Sobolev function and using the density 
of bounded functions as a first step. 
However, in the super-dimensional case $p>n$, the classical Morrey inequality 
implies that every Sobolev function in 
$W^{1,p}$ is Hölder continuous with exponent $1-\frac np$. Thus the truncation argument 
gives no improved regularity. We propose a new tool with similar approximating 
properties to the truncation for the super-dimensional setting. 
Specifically, we show that $W^{1,p}$-functions can be uniformly approximated 
by Hölder continuous functions with improved exponent $\frac {p-n}{p-1}>1-\frac np$. 
In the one-dimensional case, which is also included, the approximating function is Lipschitz for every $p>1$.

\begin{lemma}[Improved Hölder approximation]\label{lem:key}
Let $n\ge 1$, $p\in (n, \infty)$ and $\gamma_* := \frac {p-n}{p-1}\in (0, 1]$. 
Suppose $u\in W^{1,p}(\mathbb R^n)$ has compact support. 
Then there exist $u_k\in W^{1,p}(\mathbb R^n)\cap C^{0,\gamma_*}(\mathbb R^n)$ such that 
$|u_k|\le |u|$, $\nabla u_k = \chi_{\R^n\setminus A_k}\nabla u$, where $|A_k|\to 0$, and $\|u_k-u\|_{L^\infty(\R^n)} \to 0$.
Moreover, if $u$ is the zero extension of a function in $W^{1,p}_0(\Omega)$, the functions $u_k$ can be chosen with $\supp u_k\Subset \Omega$.
\end{lemma}

The main ingredient of Lemma~\ref{lem:key} is \emph{flattening}. 
Following ideas of Giacomini and Trebeschi \cite[Lemmas~3.2 and 3.3]{GiaT07},
we excise a small set of
values from the range of $u$, rather than changing $u$ directly in terms of its domain.
Specifically, the set where the $p$-energy concentrates has a small image under $u$. 
A Lipschitz change of values flattens this image, producing $C^{0,\gamma_*}$ approximations whose gradients equal either $\nabla u$ or zero. 
The coarea formula shows that the removed energy tends to zero, even without higher integrability of $\nabla u$.
The image estimate also follows from the Hausdorff-content distortion estimate of Cianchi, Korobkov and Kristensen~\cite{CiaKK25}.

We will apply this lemma to prove various results about non-standard growth energies. 
Our main application is a sharp improvement of smooth approximation in the super-dimensional case 
which  excludes the Lavrentiev gap up to the known counter-example threshold for scalar double phase energies. 
We first treat integrands with $(p,q)$-growth and H\"older dependence on the spatial variable, without assuming either convexity or continuity. 
We then develop the generalized Orlicz formulation, obtaining both a general density theorem and its super dimensional refinement under local balance conditions. 
Finally, the same flattening idea improves H\"older regularity of scalar weak subsolutions and gradient regularity of double phase minimizers in the super dimensional case.

\subsection*{Statement of results}
Let $\Omega\subset\R^n$ be open, $n\ge2$, and $F\in\Phiw(\Omega)$, where $\Phiw(\Omega)$ is the class of generalized weak $\Phi$-functions from Definition~\ref{def:weakPhi}.
We assume the \emph{$(p, q)$-growth} \cite{Mar89,Mar91,Mar93} and $\alpha$-Hölder continuity conditions
\begin{equation}\label{eq:pq}
\begin{cases}
t^p-1\lesssim F(x,t)\lesssim t^q+1,\\
|F(x,t)-F(y,t)|\lesssim |x-y|^\alpha(1+t^q), 
\end{cases}
\qquad 1\le p\le q.
\end{equation}
The implicit constants are independent of $x,y\in\Omega$ and $t\ge0$.
As usual, $f\lesssim g$ means that $f\le c g$ with constant $c\ge 1$ independent of the approximation parameters, and $f\approx g$ means that both $f\lesssim g$ and $g\lesssim f$ hold.
We further assume \adec{q} (see Definition~\ref{def:almost}), which implies, with its constant $L\ge1$, the following doubling condition
\[
 F(x,2t)\le 2^qL  F(x,t),\qquad t\ge0.
\]
However, we do not impose \ainc{p} on $F$. In particular, we do not assume convexity or continuity 
of $F(x,\cdot)$.
We then write
\[
 \mathcal F(u,\Omega)= \mathcal F_F(u,\Omega):=\int_\Omega F(x,|\nabla u|)\,dx.
\]
The main example is the \emph{double phase} function
\[
 H(x,t)=t^p+a(x)t^q,
 \qquad 1\le p\le q
 \ \ \text{and}\ \ 
 0\le a\in C^{0,\alpha}(\overline\Omega).
\]
If $\Omega$ is unbounded, the condition $a\in C^{0,\alpha}(\overline\Omega)$ means that $a$ is bounded and has a finite global H\"older seminorm. 
The growth of $H$ equals $p$ or $q$ according to whether $a$ vanishes or is positive. 
In particular, $H$ satisfies \inc{p} and \dec{q} and is doubling.
The range of exponents considered in the main theorem is
\begin{equation}\label{eq:assumptions}
 \alpha\in(0,1]
 \quad\text{and}\quad
 n<p\le q\le p+\alpha\frac{p-1}{n-1}.
\end{equation}

We derive smooth approximations with convolutions, as usual. The convergence of these convolutions  in $(p, q)$-growth functionals depends on the size of the gap $q-p$. 
The H\"older exponent
\begin{equation}\label{def:gamma*}
 \gamma_*:=\frac{p-n}{p-1}
 =1-\frac{n-1}{p-1}\in(0,1)
\end{equation}
will play a central role. It is larger than the exponent $1-\frac np$ given by the usual Sobolev embedding. 
We do not assert that every finite energy function has this better regularity. 
Rather, by Lemma~\ref{lem:key} such functions can be approximated by functions with this regularity. 
Since we have better control of the functional, we are able to consider a larger gap $q-p$. 
Specifically, for a $\gamma$-H\"older continuous function, the remaining mollification estimate requires $(1-\gamma)(q-p)\le\alpha$. 
Taking $\gamma=\gamma_*$ gives the range \eqref{eq:assumptions}.

\begin{theorem}[Smooth approximation in the super-dimensional case] \label{thm:pq}
Assume that $F\in\Phiw(\Rn)$ satisfies \adec{q} and $(p, q)$-growth \eqref{eq:pq} 
with Hölder exponent as in \eqref{eq:assumptions}.
Then, for every $u\in W^{1,p}(\R^n)$ with $\mathcal{F}(u,\R^n)<\infty$, there exists a sequence
$(v_k)\subset C_0^\infty(\Rn)$ such that
\[
\|u-v_k\|_{W^{1,p}(\R^n)} +\int_{\R^n}F(x,|\nabla u-\nabla v_k|)\,dx\to0,
\]
as $k\to \infty$. 
\end{theorem}

\begin{remark}\label{rem:extension}
When $q=p$, Hölder continuity in \eqref{eq:pq} is not needed in the previous result.
Moreover, if $u$ is the zero extension of a function in $W^{1,p}_0(\Omega)$, the functions $v_k$ can be chosen with $\supp v_k\Subset \Omega$.
\end{remark}

For prescribed boundary values, we consider the Lavrentiev gap.
For $q>p$ and fixed boundary data $g$, a Lavrentiev gap occurs if
\[
 \inf_{v\in g+W^{1,p}_0(\Omega)}\mathcal F(v,\Omega)
 <\inf_{v\in g+W^{1,q}_0(\Omega)}\mathcal F(v,\Omega).
\]
Thus restricting the admissible class changes the infimum of the energy, even though the boundary data are the same. 
In particular, some finite energy competitors cannot be approximated in energy by smoother competitors.
Zhikov studied this phenomenon in homogenization and nonstandard-growth problems~\cite{Zhi86,Zhi95}. 
Esposito, Leonetti and Mingione~\cite{ELM04} extended his two-dimensional double phase example to H\"older coefficients in arbitrary dimensions, obtaining gaps for $1<p<n<n+\alpha<q$ and connecting energy approximation with regularity. 
Fonseca, Mal\'y and Mingione~\cite{FMM04} constructed scalar double phase minimizers with fractal singular sets. 
Baroni, Colombo and Mingione~\cite{BarCM18} subsequently established regularity and energy approximation in terms of an a priori H\"older exponent of the minimizer.
Balci, Diening and Surnachev~\cite{BalDS} introduced fractal constructions for non-density and Lavrentiev gaps without requiring the exponents to lie on opposite sides of the dimension. 
For double phase growth, their Theorem~34 gives counterexamples to the density of smooth functions 
when
\[
 q>p+\alpha\max\Big\{1,\frac{p-1}{n-1}\Big\}.
\]
On the positive side, Bul\'i\v cek, Gwiazda and Skrzeczkowski~\cite{BulGS} proved smooth approximation and absence of the Lavrentiev gap under
\[
 q\le p+\alpha\max\Big\{1,\frac pn\Big\}.
\]
Their argument first reduces the problem to bounded functions via truncation 
and then uses mollification estimates. 
These two bounds agree for $p\le n$, but differ for $p>n$ when truncation loses its power. 
Corollary~\ref{cor:dp} closes the gap between these results by showing that density also holds in 
the range
\[
 p+\alpha\frac pn<q\le p+\alpha\frac{p-1}{n-1},\qquad p>n.
\]
Now the upper bound coincides with the bound from the examples from \cite{BalDS}, 
which shows the optimality of the upper bound on $q$ for density with arbitrary nonnegative $C^{0,\alpha}$ coefficients in the scalar double phase class. 

Since our construction applies also to more general functions of $(p, q)$-growth, we 
first state a more general result. 
The construction in the proof of Theorem~\ref{thm:pq} preserves prescribed boundary values and leads to the absence of the Lavrentiev gap. 
It uses the zero extension of $u-g\in W^{1,p}_0(\Omega)$, so no regularity of $\partial\Omega$ is required.

\begin{theorem}[Absence of the Lavrentiev phenomenon]\label{cor:main}
Let $\Omega\subset\Rn$ be open and bounded. Assume that $F\in\Phiw(\Rn)$ satisfies \adec{q} 
and $(p, q)$-growth \eqref{eq:pq} 
with Hölder exponent as in \eqref{eq:assumptions}, and let $g\in W^{1,q}(\Omega)$. 
Then
\[
 \inf_{v\in g+W^{1,p}_0(\Omega)}\mathcal F(v,\Omega)
 =\inf_{v\in g+W^{1,q}_0(\Omega)}\mathcal F(v,\Omega)
 =\inf_{v\in g+C_0^\infty(\Omega)}\mathcal F(v,\Omega).
\]
\end{theorem}

Applying these two theorems to the double phase function $H$ gives the following consequence.

\begin{corollary}[Density for double phase growth]\label{cor:dp}
Consider the double phase function $H$ with Hölder exponent as in \eqref{eq:assumptions}.
Then, for every $u\in W^{1,p}(\R^n)$ with $\int_{\R^n}H(x,|\nabla u|)\,dx <\infty$, there exists a sequence $(v_k)\subset C_0^\infty(\Rn)$ such that
\[
 \|u-v_k\|_{W^{1,p}(\R^n)}+\int_{\R^n}H(x,|\nabla u-\nabla v_k|)\,dx\to0,
\]
as $k\to \infty$. 
Moreover, if $\Omega\subset\Rn$ is open and bounded and $g\in W^{1,q}(\Omega)$, then
\[
 \inf_{v\in g+W^{1,p}_0(\Omega)}\mathcal F_{H}(v,\Omega)
 =\inf_{v\in g+W^{1,q}_0(\Omega)}\mathcal F_{H}(v,\Omega)
 =\inf_{v\in g+C_0^\infty(\Omega)}\mathcal F_{H}(v,\Omega).
\]

\end{corollary}

Different coefficient assumptions allow other ranges.
Borowski, Chlebicka, De Filippis and Miasojedow~\cite{BorCDM24} used the decay of weights near their zero sets, while Borowski~\cite{Bor25} treated smoother weights $a\in C^{k,\alpha}$ under $q\le p+(k+\alpha)\max\{1,\frac pn\}$.
Without H\"older continuity of $a$, Adamadze, Diening, Kopaliani and the second author~\cite{AdaDKO26} obtained smooth density and H\"older regularity of solutions under a Muckenhoupt-type condition. 
Their condition includes bounded H\"older coefficients with $\frac qp\le 1+\frac\alpha n$ and also permits discontinuous or unbounded coefficients.

Generalized Orlicz spaces encompass Orlicz, variable exponent and double phase growth; see \cite{HarH19}. In this setting the inequalities between the global growth rates $p$ and $q$ are 
replaced by a local balance or comparison condition. 
 A basic density result is given in \cite[Theorem~6.4.4]{HarH19}. 
Modular approximation under balance conditions was developed by Ahmida, Chlebicka, Gwiazda and Youssfi~\cite{AhmCGY18}, and extended to anisotropic integrands without doubling by Borowski and Chlebicka~\cite{BorC22};
see also~\cite{BorCM24} for Lavrentiev gaps.
Bul\'i\v cek, Gwiazda and Skrzeczkowski~\cite{BulGS} proved smooth density and absence of the Lavrentiev gap in the Musielak--Orlicz setting under a local comparison condition and suitable growth assumptions.
The generalized Jensen inequalities in~\cite{Has23} relate local balance conditions to convex minorants.
For manifold-valued maps, Antonini, De Filippis and Pacchiano Camacho~\cite{AntDPC26} recently studied smooth approximation and Lavrentiev gaps, including the role of topological constraints. 
We expect that these results can also be improved in the super-dimensional case 
by our method.

We establish generalized Orlicz density results in Theorem~\ref{thm:orlicz-super}.
This result is new also in the sub-dimensional case.
Borowski and Chlebicka~\cite{BorC22} allow fully anisotropic and non-doubling growth, but impose a stronger balance condition of \aone{}-type based on \cite{Has23}.
(To exemplify the assumptions, for double phase growth \aone{} corresponds to $q-p\le\alpha \frac pn$, whereas \aonen{n} allows $q-p\le\alpha$.)
The range $q-p\le\alpha$ is covered by \cite{BulGS}, but for the generalized  Orlicz case they resort to the stronger balance condition; our first contribution is  a generalization of \cite{BulGS} to the generalized Orlicz case.
In the scalar super-dimensional case, when \ainc{p} holds with $p>n$, 
Theorem~\ref{thm:orlicz-super} requires only the even weaker condition \aonen{s_*}, where
$s_*:=\frac{n(p-1)}{n-1}>n$.
For double phase growth, this yields $q-p\le\alpha \frac{p-1}{n-1}$, improving the bound $q-p\le\alpha \frac pn$ in~\cite{BulGS} in the super-dimensional case.
\comment{Peter: \cite{BulGS} deals with Musielak--Orlicz--Sobolev space only under \aonen{p}.}

Related balance conditions arise in elliptic regularity. 
Colombo and Mingione~\cite{ColM15} developed double phase regularity theory, while Baroni, Colombo and Mingione~\cite{BarCM16} studied borderline and logarithmic growth. 
De Filippis and Mingione studied regularity under energy approximation assumptions~\cite{DeFM20} and obtained improved growth bounds for H\"older continuous minimizers~\cite{DeFM21}.
See also~\cite{DeFM23,DeFM25} for sharp Schauder estimates in the nonuniformly elliptic setting.
For generalized Orlicz growth, we obtained regularity under spatial balance conditions~\cite{HasO22,HasO22b} and examined their relaxation using a priori information~\cite{HasO23}.
Density permits energy-space test functions to be approximated by smooth ones.
In Section~\ref{sect:maximal}, our flattening argument also gives $C^{0,\gamma_*}_{\loc}$-regularity for scalar weak subsolutions.
For local double phase minimizers, combining this estimate with \cite[Theorem~1.3]{BarCM18} yields maximal regularity for the improved range
\[
 n<p<q<p+\alpha\frac{p-1}{n-1}.
\]
The gradient regularity conclusion uses a strict upper bound, while the density theorem includes the endpoint.

This paper is organized as follows. 
In Section~\ref{sect:approximation}, we prove the flattening result, Lemma~\ref{lem:key}. 
Section~\ref{sect:superdimensional} proves the super-dimensional approximation and
Lavrentiev results. Section~\ref{sect:orlicz} treats generalized Orlicz space density, 
and Section~\ref{sect:maximal} gives the maximal regularity application.

%%%%%%%%%%%%%%%%%%%%%%%%%%%%%%%%%%%%%%%%%%%%%%%%%%%%%%%%%%%%%%%%%%%%%%%%%%%%%%%%%%%
%%%%%%%%%%%%%%%%%%%%%%%%%%%%%%%%%%%%%%%%%%%%%%%%%%%%%%%%%%%%%%%%%%%%%%%%%%%%%%%%%%%
%%%%%%%%%%%%%%%%%%%%%%%%%%%%%%%%%%%%%%%%%%%%%%%%%%%%%%%%%%%%%%%%%%%%%%%%%%%%%%%%%%%
\section{Proof of the approximation lemma}\label{sect:approximation}

We write $B_r(x)$ for the open ball with center $x$ and radius $r$, and omit the center when it is fixed. 
For a ball $B$, $\lambda B$ denotes the concentric ball with radius multiplied by $\lambda$.
We use $u_B:=\fint_B u\,dx=|B|^{-1}\int_B u\,dx$ for the average. 
In this section we consider also the one-dimensional case and only assume that $n\ge 1$.

We begin with two forms of Morrey's inequality. 
The following form can be found for $C^1$-functions, e.g., in \cite[Theorem~5.29]{Kin26}; the statement for continuous Sobolev functions follows by approximation. 
If $p>n$ and $u\in W^{1,p}(B_r)\cap C(B_r)$, then 
\begin{equation}\label{eq:morrey-oscillation}
\osc_{B_r}\,u := \sup_{B_r}u-  \inf_{B_r} u
 \le c r^{1-n/p}\left(\int_{B_r}|\nabla u|^p\,dx\right)^{1/p}.
\end{equation}
In particular, for $u\in W^{1,p}(\Rn)\cap C(\Rn)$ and $x\ne z$, this inequality 
on $B_{2|x-z|}(z)$ gives
\begin{equation}\label{eq:morrey-twopoints}
 |u(x)-u(z)|
 \le  c \,|x-z|^{1-n/p}
 \left(\int_{B_{2|x-z|}(z)}|\nabla u|^p\,dy\right)^{1/p}.
\end{equation}
Here and below, $c$ may depend on $n$ and $p$. 
The larger ball in the second estimate is only used to include both points.

We also recall the scalar Sobolev coarea formula \cite[Proposition~2.1]{BZ}.
Let $u\in W^{1,1}(\R^n)$, represented by its precise Sobolev
representative, with arbitrary values on its $\mathcal H^{n-1}$-null
exceptional set.
For every open $\Omega\subset\R^n$ and every nonnegative function $f$ on $\R^n$, we have
\begin{equation}\label{eq:weighted-sobolev-coarea}
 \int_\Omega f(x)\,|\nabla u(x)|\,dx
 =
 \int_{\mathbb R} \int_{u^{-1}(t)} f(x)\,d\mathcal H^{n-1}(x)\, dt,
\end{equation}
where $u^{-1}(t)=\{x\in\Omega:u(x)=t\}$; here both sides are understood as nonnegative integrals,
possibly infinite, with the convention $0\cdot\infty=0$.

The following lemma is a standard consequence of this formula.
A closely related result characterizing sets $A$ with 
$\big|u^{-1}(A)\cap\{|\nabla u|>0\}\big|=0$ 
appears in \cite[Corollary~4.6]{LeoniMorini07}.
We include a direct proof for completeness. 
Note that the choice of $\delta$ depends on both $\epsilon$ and $u$.

\begin{lemma}
\label{lem:inversset}
Let $n\ge 1$ and let $u\in W^{1,1}(\mathbb R^n)$
have compact support. For every $\epsilon>0$ there exists $\delta>0$ such that
\[
\big|u^{-1}(A)\cap\{|\nabla u|>0\} \big|<\epsilon,
\]
whenever $A\subset\mathbb R$ is Lebesgue measurable and satisfies
$\mathcal L^1(A) <\delta$.
\end{lemma}
\begin{proof}
Set $G:=\{x\in\mathbb R^n:|\nabla u(x)|>0\}$.
Since $u$ has compact support, $|G|<\infty$. Define 
\[
 g(x):=
 \begin{cases}
  |\nabla u(x)|^{-1}, & x\in G,\\
  0,                 & x\notin G.
 \end{cases}
\]
Applying \eqref{eq:weighted-sobolev-coarea} with $f=g$
and $\Omega=\mathbb R^n$, we obtain
\[
 \int_{\mathbb R}   \int_{u^{-1}(t)\cap G}
      \frac{1}{|\nabla u(x)|}\,d\mathcal H^{n-1}(x)\, dt
 = \int_{\R^n}g\,|\nabla u|\, dx=  |G|<\infty.
\]
Consequently, the function
\[
 h(t):=
 \int_{u^{-1}(t)\cap G} \frac{1}{|\nabla u(x)|}\,d\mathcal H^{n-1}(x)
\]
belongs to $L^1(\mathbb R)$. For any Borel set $A\subset\mathbb R$, applying 
\eqref{eq:weighted-sobolev-coarea} with $f=g\chi_{u^{-1}(A)}$ gives
\[
\big|u^{-1}(A)\cap G\big|
=
\int_{\Rn} \chi_{u^{-1}(A)} \,g\,|\nabla u|\, dx
=
\int_A h(t)\,dt 
\le 
K\mathcal L^1(A)+\int_{\{t:h(t)>K\}}h(t)\,dt,
\]
for every $K>0$. Since $h\in L^1(\R)$, there exists large $K>0$ such that
 $\int_{\{t:h(t)>K\}}h(t)\,dt<\frac{\epsilon}{2}$. Therefore, by  choosing 
 $\delta:=\frac{\epsilon}{2K}$, we obtain the claim.
\end{proof}

We now prove the key lemma. 

\begin{proof}[Proof of Lemma~\ref{lem:key}]
We divide the proof into four steps.

\smallskip\noindent\textbf{Step 1: Concentration set and its Vitali covering.}
If $\|\nabla u\|_{L^p(\mathbb R^n)}=0$, continuity and compact support imply $u\equiv 0$. In this case, we can simply choose $u_k\equiv 0$ as an approximating sequence. Assume instead that $A:=\|\nabla u\|_{L^p(\mathbb R^n)}^p>0$.
For $k\ge 1$, we define
\[
 E_k:=\left\{x\in\mathbb R^n:
      \sup_{r>0}r^{-\gamma_*}\int_{B_r(x)}|\nabla u|^p\, dy>k\right\}.
\]
Since the integral over $B_r(x)$ is continuous in $x$, the set $E_k$ is open.

Let first $z\in \Rn\setminus E_k$. Then 
\[
 \int_{B_r(z)}|\nabla u|^p\, dx\le k r^{\gamma_*}\quad\text{for every }\ r>0.
\]
For any $x\ne z$, Morrey's inequality \eqref{eq:morrey-twopoints} yields
\begin{equation}\label{eq:good-anchor}
 \begin{split}
 |u(x)-u(z)|
 &\le c \,|x-z|^{1-n/p}\left(\int_{B_{2|x-z|}(z)}|\nabla u|^p\, dx\right)^{1/p}\\
 &\le c k^{1/p}|x-z|^{1-n/p+\gamma_*/p} = c k^{1/p}|x-z|^{\gamma_*};
 \end{split}
\end{equation}
in the last equality, we used the relation $1-n/p+\gamma_*/p=\gamma_*$ from \eqref{def:gamma*}. 
We conclude that $u$ is locally $\gamma_*$-H\"older continuous at every $z\in \Rn\setminus E_k$.

We next consider $z\in E_k$ and observe that 
\[
 R_z:=\sup\left\{r>0:\int_{B_r(z)}|\nabla u|^p \,dx>r^{\gamma_*} k \right\}\le\Big(\frac A k\Big)^{1/\gamma_*}<\infty.
\]
By the definition of supremum, we can choose a radius $r_z\in (\frac12 R_z,R_z]$ such that 
\begin{equation}\label{eq:violating-radius}
 \int_{B_{r_z}(z)}|\nabla u|^p\, dx>r_z^{\gamma_*} k;
\end{equation}
since $2r_z > R_z$, we also have that 
\begin{equation}\label{eq:large-good-radius}
 \int_{B_r(z)}|\nabla u|^p\,dx\le r^{\gamma_*} k \quad \text{for all }\ r\ge 2r_z.
\end{equation}
As the radii $r_z$ are uniformly bounded, Vitali's covering lemma ensures the existence of a sequence $(z_i)_{i=1}^\infty \subset E_k$ such that the balls $B_i:=B_{r_i}(z_i)$ are mutually disjoint and $E_k\subset\bigcup_i 5B_i$, where $r_i:=r_{z_i}$. 
By \eqref{eq:violating-radius} and disjointness,
\begin{equation}\label{eq:radii-sum}
 \sum_i r_i^{\gamma_*}\le\frac{1}{k}\sum_i\int_{B_i}|\nabla u|^p\, dz
 \le\frac{A}{k}.
\end{equation}
On the other hand, using Morrey's inequality \eqref{eq:morrey-oscillation} and
\eqref{eq:large-good-radius} we obtain
\[
 \osc_{5B_i} u
 \le c r_i^{1-n/p}\left(\int_{5B_i}|\nabla u|^p\, dz\right)^{1/p}
 \le c k^{1/p}r_i^{(p-n+\gamma_*)/p}.
\]
In the next step, we want to estimate the right-hand side using \eqref{eq:radii-sum}; 
for this we need that the exponents of $r_i$ satisfy $(p-n+\gamma_*)/p=\gamma_*$, which is exactly \eqref{def:gamma*}.

\smallskip\noindent\textbf{Step 2: Bounding the image of the concentration set.}
It is not guaranteed that $u(E_k)$ is open. Thus, we construct a slightly larger open set containing $(-\frac1k,\frac1k) \cup u(E_k)$.
Define 
\[
 I_k :=\big(-\tfrac1k,\tfrac1k\big)\,
       \cup\, \bigcup_i \Big(\inf_{5B_i}u-\tfrac1{2^ik},\sup_{5B_i}u+\tfrac1{2^ik}\Big).
\]
Note that $I_k$ is open. Since $E_k \subset \bigcup_i 5B_i$, we have 
\[
\big(-\tfrac1k,\tfrac1k\big) \cup u(E_k)  \, \subset \, I_k 
\]
and by \eqref{eq:radii-sum}, 
\begin{equation}\label{eq:I2}
\mathcal{L}^1(I_k) \le \frac{2}{k}+ \sum_i\big(\osc_{5B_i}u + \tfrac2{2^ik}\big)
 \le \frac{4}{k}+c k^{1/p}\sum_i r_i^{\gamma_*} 
 \le\frac{4}{k} +  c \frac{ A}{k^{1/p'}}.
\end{equation}

\smallskip\noindent\textbf{Step 3: Construction of the H\"older continuous approximation.}
Define
\[
 u_k(x):=T_k(u(x)), 
 \quad\text{where }\ T_k(t):=\int_0^t\chi_{\R\setminus I_k}(s)\,ds.
\]
Note that $T_k$ is $1$-Lipschitz, $|T_k(t)|\le |t|$, 
\[
t-T_k(t)= \int_0^t\chi_{I_k}(s)\,ds,\quad\text{and}\quad T_k(t)=0\ \ \text{for}\ |t|\le \tfrac1k.
\]
In particular, $T_k$ is constant on every interval contained in $I_k$.

If $\mathcal C$ is a connected component of the open set $E_k$, the continuity of $u$ implies that $u(\mathcal C)$ is an interval contained in $I_k$. Thus, $u_k=T_k \circ u$ is constant on $\mathcal C$.
Consider distinct $x,y\in\mathbb R^n$.
If the segment $[x,y]$ lies entirely within a connected component of $E_k$, then $u_k(x)-u_k(y)=0$. Otherwise, the segment contains a point $z\in \Rn\setminus E_k$.
The Lipschitz property of $T_k$ combined with the Hölder estimate \eqref{eq:good-anchor} gives
\[
 \begin{split}
 |u_k(x)-u_k(y)| 
 \le |u(x)-u(y)| 
 \le  |u(x)-u(z)|+|u(z)-u(y)|
 &\lesssim k^{1/p}\left(|x-z|^{\gamma_*}+|z-y|^{\gamma_*}\right)\\
 &\lesssim k^{1/p} |x-y|^{\gamma_*}.
 \end{split}
\]
Hence, $u_k\in C^{0,\gamma_*}(\mathbb R^n)$.

\smallskip\noindent\textbf{Step 4: Convergence.}
From $|T_k(t)|\le |t|$ it follows that $|u_k|\le |u|$.
By the definition of $T_k$ and \eqref{eq:I2},
\[
|u(x)-u_k(x)| 
= 
|u(x)-T_k(u(x))|
=
\bigg|\int_0^{u(x)} \chi_{I_k}(s)\,ds\bigg|
\le 
\mathcal{L}^1(I_k)
\le 
\frac{4}{k} +  c \frac{ A}{k^{1/p'}}.
\]
The right-hand side can be made arbitrarily small by taking $k$ large enough. 
Hence $u_k\to u$ uniformly. 

Since $T_k$ is Lipschitz, the chain rule yields
\[
 \nabla u_k =\chi_{\mathbb R \setminus I_k}(u) \nabla u = \chi_{\{u\notin I_k\}} \nabla u
 \quad\text{and}\quad
\nabla u  -   \nabla u_k = \chi_{u^{-1}( I_k)} \nabla u
\]
almost everywhere. Define $A_k:=u^{-1}(I_k)\cap\{|\nabla u|>0\}$.
By Lemma~\ref{lem:inversset}, the measure $|A_k|$ is  small if $\mathcal{L}^1(I_k)$ is small. In view of \eqref{eq:I2}, this is achieved by choosing $k$ 
large.

Finally, suppose $w\in W^{1,p}_0(\Omega)$ and $u$ is its zero extension.
Then $u=0$ in $\Rn\setminus\Omega$ except in a set of capacity zero \cite[Theorem~4.33]{Kin26}, 
but since $p>n$ the exceptional set must be empty. 
Since $T_k(t)=0$ for $|t|\le \frac1k$ and $u$ is continuous, we deduce that
\[
 \supp u_k \subset\big\{x\in\overline\Omega:|u(x)|\ge \tfrac1k\big\}
 \Subset\Omega. \qedhere
\]
\end{proof}

%%%%%%%%%%%%%%%%%%%%%%%%%%%%%%%%%%%%%%%%%%%%%%%%%%%%%%%%%%%%%%%%%%%%%%%%%%%%%%%%%%%
%%%%%%%%%%%%%%%%%%%%%%%%%%%%%%%%%%%%%%%%%%%%%%%%%%%%%%%%%%%%%%%%%%%%%%%%%%%%%%%%%%%
%%%%%%%%%%%%%%%%%%%%%%%%%%%%%%%%%%%%%%%%%%%%%%%%%%%%%%%%%%%%%%%%%%%%%%%%%%%%%%%%%%%
\section{Density in the super-dimensional case}\label{sect:superdimensional}

In this section, we show that higher H\"older continuity of Sobolev functions 
allows us to mollify under a weaker
balance condition in line with prior results \cite{BarCM18, BenHHK21,HasO23}. 
We start with introducing some conventions and definitions. 
From here onward we assume that $n\ge 2$.

We fix a standard mollifier
$\eta\in C_0^\infty(B_1(0))$ with $\eta\ge0$ and $\int_\Rn\eta\,dx=1$.
For $r>0$, set $\eta_r(x):=r^{-n}\eta(\frac xr)$.
We denote cut-off functions by $\zeta$.
We define growth conditions and generalized
weak $\Phi$-functions following \cite{HarH19}.

\begin{definition}\label{def:almost}
Let $f: \Omega\times [0,\infty) \to [0,\infty]$ and $p,q>0$. We say that $f$ satisfies 
\begin{itemize}
\item[(aInc)$_p$]\label{ainc}\label{inc}
if there exists $L\ge 1$ such that $f(x,\lambda t)\le L \lambda^p f(x,t)$ 
for all $t\ge 0$, $0<\lambda \le 1$ and a.e.\ $x\in \Omega$;
\item[(aDec)$_q$]\label{adec}\label{dec}
if there exists $L\ge 1$ such that $f(x,\Lambda t)\le L \Lambda^q f(x,t)$ 
for all $t\ge 0$, $\Lambda \ge 1$ and a.e.\ $x\in \Omega$.
\end{itemize}
If the condition holds with $L=1$, then we use \inc{p} and \dec{q}. 
\end{definition}

\begin{definition}\label{def:weakPhi}
We say that $\phi:\Omega\times[0,\infty)\to[0,\infty]$ is a
\textit{generalized weak $\Phi$-function}, denoted by $\phi\in\Phiw(\Omega)$,
if $x\mapsto\phi(x,|f(x)|)$ is measurable for every measurable $f$ on
$\Omega$, $\phi(x,\cdot)$ is nondecreasing and satisfies \ainc{1}
uniformly in $x$, and, for almost every $x\in\Omega$,
$\phi(x,0)=\lim_{t\to0^+}\phi(x,t)=0$ and $\lim_{t\to\infty}\phi(x,t)=\infty$.
\end{definition}

\begin{lemma}\label{lem:mollification}
Let $\Omega\subset\Rn$ be bounded and open, $1<p\le q<\infty$,
$0<\alpha\le1$, and $0<\gamma<1$.
Assume that $F\in\Phiw(\Omega)$ satisfies \adec{q} as well as the 
$(p, q)$-growth and $\alpha$-Hölder conditions from \eqref{eq:pq}.
Let $u\in W^{1,p}(\mathbb R^n)\cap C^{0,\gamma}(\mathbb R^n)$ satisfy
$K:=\supp u\Subset\Omega$ and
$\int_\Omega F(x,|\nabla u|)\,dx<\infty$.
If
\begin{equation}\label{eq:holder-balance}
 \alpha-(1-\gamma)(q-p)\ge0,
\end{equation}
then the mollifications $u_\epsilon=\eta_\epsilon*u$ satisfy
\[
 \|u_\epsilon-u\|_{W^{1,p}(\mathbb R^n)}+ \int_\Omega F(x,|\nabla u_\epsilon-\nabla u|)\,dx\longrightarrow0
 \quad\text{as }\epsilon\to0.
\]
\end{lemma}
\begin{proof}
We restrict $\epsilon$ such that
\[
 0<\epsilon<\min\{1,\tfrac13\dist(K,\partial\Omega)\}.
\]
Then $u_\epsilon\in C_0^\infty(\Omega)$,
$u_\epsilon\to u$ in $W^{1,p}(\Omega)$ and
$\nabla u_\epsilon\to\nabla u$ almost everywhere.
If $q=p$, then $\nabla u_\epsilon\to\nabla u$ in
$L^1(\Omega)\cap L^q(\Omega)$. Since \ainc{1} and the upper growth
bound in \eqref{eq:pq} yield $F(x,t)\lesssim t+t^q$, the conclusion follows directly.

From now on, we assume that $q>p$.
Continuity of the weak $\Phi$-function at $t=0$ implies that
\[
F(\cdot, |\nabla u -\nabla u_\epsilon|) \to 0
\quad  \text{a.e. as }\ \epsilon\to 0.
\]
We will show that $F(\cdot,|\nabla u -\nabla u_\epsilon|)$ is uniformly integrable. Then, Vitali's convergence theorem \cite[p.~133]{Rud87} yields the claim.

Since $\int_{\mathbb R^n}\nabla\eta_\epsilon(y)\,dy=0$,
\begin{equation}\label{eq:Dvepsilon_bound}
\begin{split}
|\nabla u_\epsilon(x)|
& =\left|\int_{\mathbb R^n}\nabla\eta_\epsilon(y)
(u(x-y)-u(x))\,dy\right| \\
&\le \epsilon^{-n-1}[u]_{C^{0,\gamma}}
\int_{\mathbb R^n}|\nabla\eta(\tfrac y\epsilon)|\,|y|^\gamma\,dy
=\epsilon^{\gamma-1}[u]_{C^{0,\gamma}}
\underbrace{\int_{\mathbb R^n}|\nabla\eta(z)|\,|z|^\gamma\,dz}_{=:c_\eta}.
\end{split}
\end{equation}
Fix $x\in\Omega$. If $\dist(x,K)>\epsilon$,
we trivially have $\nabla u(x)=\nabla u_\epsilon(x)=0$.
Thus, we may assume $\dist(x,K)\le\epsilon$,
which implies $\overline B_\epsilon(x)\subset\Omega$.
For such $x$, we define
\[
 F^-_{x,\epsilon}(t):=
 \inf_{y\in\overline B_\epsilon(x)}F(y,t)%,
\quad\text{and}\quad
 \tilde F_{x,\epsilon}:=(F^-_{x,\epsilon})^{**},
\]
where $f^{**}=(f^*)^*$ and
$f^*(s):=\sup_{t\ge0}\{st-f(t)\}$. Then, $F^-_{x,\epsilon}$ satisfies \ainc{1} and \adec{q} with the same constants as $F$, and moreover the lower bound in \eqref{eq:pq} ensures
$F^-_{x,\epsilon}\in\Phiw$.
By \cite[Lemma~2.4.2 and Proposition~2.4.5]{HarH19} and \adec{q},
$\tilde F_{x,\epsilon}$ is convex and nondecreasing, and
\[
 \tilde F_{x,\epsilon}(t)
 \le F^-_{x,\epsilon}(t)
 \lesssim \tilde F_{x,\epsilon}(t),
\]
with constants independent of $x$ and $\epsilon$.
Hence, by the $\alpha$-Hölder continuity \eqref{eq:pq},
\[
 F(x,t)
 \lesssim F^-_{x,\epsilon}(t)+\epsilon^\alpha(1+t^q)
 \lesssim \tilde F_{x,\epsilon}(t)
          +\epsilon^\alpha(1+t^q).
\]
By this estimate, Jensen's inequality for the convex functions 
$\tilde F_{x,\epsilon}$ and $|\cdot|^p$,
\eqref{eq:Dvepsilon_bound}, the $(p, q)$-growth \eqref{eq:pq}, and $|\nabla u_\epsilon|\lesssim 
\epsilon^{\gamma-1}$, we obtain
\begin{equation*}
\begin{split}
F(x,|\nabla u_\epsilon(x)|)
&\lesssim
\tilde F_{x,\epsilon}(|\nabla u_\epsilon(x)|)
+\epsilon^\alpha(1+|\nabla u_\epsilon(x)|^q)\\
&\lesssim
\big(\eta_\epsilon*
 \tilde F_{x,\epsilon}(|\nabla u|)\big)(x) +\epsilon^\alpha
\big(1+(c_\eta [u]_{C^{0,\gamma}}\epsilon^{\gamma-1})^{q-p}
|\nabla u_\epsilon(x)|^p\big)\\
&\lesssim
\big(\eta_\epsilon*F(\cdot,|\nabla u|)\big)(x) +\epsilon^\alpha\big(1+\epsilon^{(\gamma-1)(q-p)} (\eta_\epsilon*|\nabla u|^p)(x)\big)\\
&\lesssim
\big[1+\epsilon^{\alpha + (\gamma-1)(q-p)}\big]\big[\big(\eta_\epsilon*F(\cdot,|\nabla u|)\big)(x) +1\big]\\
&\lesssim \big(\eta_\epsilon*F(\cdot,|\nabla u|)\big)(x)+1,
\end{split}
\end{equation*}
where, in the last inequality, we used \eqref{eq:holder-balance} and $0<\epsilon<1$.
Extending $F(\cdot,|\nabla u|)$ outside $\Omega$ by zero,
$F(\cdot,|\nabla u|)\in L^1(\Rn)$. Thus the mollifications of $F(\cdot,|\nabla u|)$ are uniformly
integrable. Consequently, since $\Omega$ is bounded, the previous inequality ensures that
$\{F(\cdot,|\nabla u_\epsilon|)\}_\epsilon$
is uniformly integrable on $\Omega$.
Moreover, the monotonicity and doubling of $F(x,\cdot)$ yield
\[
 F(x,|\nabla u_\epsilon-\nabla u|)
 \lesssim F(x,|\nabla u_\epsilon|)+F(x,|\nabla u|).
\]
Thus $F(\cdot,|\nabla u_\epsilon-\nabla u|)$ is uniformly
integrable on $\Omega$. Therefore, the proof is completed.
\end{proof}

Before applying the preceding lemmas, we reduce to functions
with compact support by truncating small values. In contrast to the usual 
proof where we multiply $u$ by a cut-off function, this approach does not require 
us to assume that $\int_{\R^n} F(x,|u|)\,dx <\infty$.

\begin{lemma}\label{lem:dens-compact-support-F}
Let $F\in\Phiw(\Rn)$, and $u\in W^{1,p}(\Rn)$ for some $p>n$ with $\int_{\R^n} F(x,|\nabla u|)\,dx <\infty$.
Then there exist compactly supported functions
$u_m \in W^{1,p}(\Rn)$ with $\int_{\R^n} F(x,|\nabla u_m|)\,dx <\infty$ and
$\supp u_m\subset\supp u$ such that 
\[
\|u_m-u\|_{W^{1,p}(\R^n)}+\int_{\Rn} F(x,|\nabla u_m- \nabla  u|)\,dx \to 0
\]
as $m\to \infty$.
\end{lemma}
\begin{proof}
Since $p>n$, the Sobolev--Morrey embedding gives
$u\in C_0(\Rn)$.
In particular, $u$ is bounded and
$\{x\in\Rn:|u(x)|\ge \frac1m\}$ is compact for every $m\in\N$.  For $m\in\N$, define
\[
 u_m:= 
 \begin{cases}
 u -\tfrac1m,& \text{if }\ u>  \tfrac1m\\
  0,& \text{if }\ |u| \le  \tfrac1m\\
   u +\tfrac1m,& \text{if }\ u<-  \tfrac1m\\
 \end{cases}
\]
Then $|u_m|\le |u|$ and
\[
 \supp u_m\subset\{|u|\ge \tfrac1m\}\subset\supp u,
\]
so $u_m$ is compactly supported. By \cite[Corollary~2.1.8]{Zie89},
\[
 \nabla u_m=\chi_{\{|u|>1/m\}}\nabla u
 \quad\text{a.e.}
\]
Since $F(x,0)=0$, we have
\[
 \int_\Rn F(x,|\nabla u_m|)\,dx
 \le\int_\Rn F(x,|\nabla u|)\,dx<\infty.
\]
Moreover, $\nabla u=0$ almost everywhere in $\{u=0\}$, and hence
\[
 \int_\Rn F(x,|\nabla u-\nabla u_m|)\,dx =\int_{\{|u|\le 1/m\}} F(x,|\nabla u|)\,dx
 \to \int_{\{u=0 \}}F(x,|\nabla u|)\,dx =0,
\
\]
by the dominated convergence theorem. The same argument with $F(x,t)=t^p$
gives $\nabla u_m\to\nabla u$ in $L^p(\Rn)$.
Finally, $u_m\to u$ pointwise and $|u-u_m|\le|u|$, so the dominated
convergence theorem also yields $u_m\to u$ in $L^p(\Rn)$.
\end{proof}

We now combine the flattening from Lemma~\ref{lem:key} with
mollification. The choice of $\gamma_*$ is exactly what gives the
exponent range in \eqref{eq:assumptions}.

\begin{proof}[Proof of Theorem~\ref{thm:pq}]
By Lemma~\ref{lem:dens-compact-support-F},  it is enough to consider
the case that $u$ has compact support. Suppose $\supp u\Subset\Omega$ for some
bounded open set $\Omega$. Let $k>0$.
By Lemma~\ref{lem:key} and dominated convergence, there exists
$w_k\in W^{1,p}(\Rn)\cap C^{0,\gamma_*}(\Rn)$ with compact support in
$\Omega$ such that
\[
 \|u-w_k\|_{W^{1,p}(\Rn)}
 +\int_\Rn F(x,|\nabla u-\nabla w_k|)\,dx<\frac1k.
\]
Since $(1-\gamma_*)(q-p)\le\alpha$ by \eqref{eq:assumptions},
Lemma~\ref{lem:mollification} yields $v_k\in C_0^\infty(\Omega)$ such that
\[
 \|w_k-v_k\|_{W^{1,p}(\Rn)}
 +\int_\Rn F(x,|\nabla w_k-\nabla v_k|)\,dx<\frac1k.
\]
By the doubling property from \adec{q},
\[
 \int_\Rn F(x,|\nabla u-\nabla v_k|)\,dx
 \lesssim\int_\Rn F(x,|\nabla u-\nabla w_k|)
                    +F(x,|\nabla w_k-\nabla v_k|)\,dx
\lesssim\frac1k.
\]
Similarly, $\|u-v_k\|_{W^{1,p}(\Rn)}\to 0$ as
$k\to \infty$.
\end{proof}

We conclude the section by proving the absence of a Lavrentiev gap,
Theorem~\ref{cor:main}. We note that if $F(x,\cdot)$ is continuous
for a.e.~$x$, the conclusion follows easily from the local approximation
construction in Theorem~\ref{thm:pq} and Vitali's convergence theorem \cite[p.~133]{Rud87}.
However, under our assumptions, $F(x,\cdot)$ may be discontinuous,
so a more delicate argument is needed.
To this end, we first apply a small local contraction to $u$,
preserving its boundary values, and then approximate the resulting
function. The contraction allows us to use the monotonicity of
$F(x,\cdot)$ to control the energy without assuming continuity.

\begin{proof}[Proof of Theorem~\ref{cor:main}]
Let $u\in g+W^{1,p}_0(\Omega)$ satisfy
$\mathcal F(u,\Omega)<\infty$.
Fix $\epsilon>0$ and choose a compact set $K\subset\Omega$ such that
\[
 \int_{\Omega\setminus K}F(x,|\nabla u|)\,dx<\epsilon.
\]
Let $\zeta\in C_0^\infty(\Omega)$ satisfy $0\le\zeta\le1$ and
$\zeta=1$ in a neighborhood of $K$. Since $p>n$, $u$ is bounded on
$\supp\zeta$. For $0<\delta<1$, set
\[
 u_\delta:=(1-\delta\zeta)u.
\]
Then
$\nabla u_\delta=(1-\delta\zeta)\nabla u-\delta u\nabla\zeta$
and $u_\delta\in g+W^{1,p}_0(\Omega)$.
By doubling, \ainc{1}, and the upper growth bound in \eqref{eq:pq},
\[
 F(x,|\nabla u_\delta|)
 \lesssim F(x,|\nabla u|)+\delta\big(1+\|u\|_{L^\infty(\supp\zeta)}^q \|\nabla\zeta\|_\infty^q\big).
\]
Hence
\begin{equation}\label{eq:u_delta}
 \int_\Omega F(x,|\nabla u_\delta|)\,dx<\infty
 \quad\text{and}\quad
 \int_{\Omega\setminus K}F(x,|\nabla u_\delta|)\,dx
 \le c\epsilon+c_0\delta.
\end{equation}
Here $c_0$ may depend on $u$, $\zeta$ and $\Omega$, but is independent of $\delta$.

Set $w_\delta:=u_\delta-g\in W^{1,p}_0(\Omega)$.
Since $g\in W^{1,q}(\Omega)$ and $F(x,|\nabla w_\delta|)\lesssim F(x,|\nabla u_\delta|) + |\nabla g|^q+1$, we have that $\mathcal F(w_\delta,\Omega)<\infty$.
By Theorem~\ref{thm:pq} and Remark~\ref{rem:extension} for the zero extension
of $w_\delta$, we obtain $v_j\in g+C_0^\infty(\Omega)$ for sufficiently small mollification radii such that
\begin{equation}\label{eq:v_j}
 v_j\to u_\delta \quad\text{in }W^{1,p}(\Omega),
 \qquad
 \int_\Omega F(x,|\nabla v_j-\nabla u_\delta|)\,dx\to0
\end{equation}
as $j\to \infty$. Passing to a subsequence, we also have
$\nabla v_j\to\nabla u_\delta$ almost everywhere.
By monotonicity and doubling of $F$,
\[
 F(x,|\nabla v_j|)
 \lesssim F(x,|\nabla u_\delta|)
       +F(x,|\nabla v_j-\nabla u_\delta|),
\]
so $F(\cdot,|\nabla v_j|)$ is uniformly integrable.

Define $g_j:=\left(F(\cdot,|\nabla v_j|)-F(\cdot,|\nabla u|)\right)_+$ on $K$.
Then $0\le g_j\le F(\cdot,|\nabla v_j|)$ and
$F(\cdot,|\nabla v_j|)\le F(\cdot,|\nabla u|)+g_j$.
Since $|\nabla v_j|\to(1-\delta)|\nabla u|$ a.e.\ in $K$, for a.e.\ 
$x\in K$ with $|\nabla u(x)|>0$ there exists $j_x\in\N$ such that
$|\nabla v_j(x)|\le|\nabla u(x)|$ for all $j\ge j_x$.
By monotonicity of $F$, $g_j(x)=0$ for these $j$.
At points where $\nabla u(x)=0$, we use
$\lim_{t\to0}F(x,t)=F(x,0)=0$. Thus $g_j\to0$ a.e. in $K$.
The preceding uniform integrability and Vitali's convergence theorem \cite[p.~133]{Rud87} give
\[
 \lim_{j\to\infty}\int_K g_j\,dx=0.
\]
Consequently,
\[
 \limsup_{j\to\infty}\int_K F(x,|\nabla v_j|)\,dx
 \le\int_K F(x,|\nabla u|)\,dx
       +\lim_{j\to\infty}\int_K g_j\,dx
 =\int_K F(x,|\nabla u|)\,dx.
\]

On $\Omega\setminus K$, the preceding estimates \eqref{eq:u_delta} and \eqref{eq:v_j} yield
\[
 \limsup_{j\to\infty}
 \int_{\Omega\setminus K}F(x,|\nabla v_j|)\,dx
 \le c \epsilon+c_0\delta.
\]
Therefore, since $v_j\in g+C_0^\infty(\Omega)$, it follows that
\[
 \inf_{v\in g+C_0^\infty(\Omega)}\mathcal F(v,\Omega)
 \le \mathcal F(u,\Omega)+c\epsilon+c_0\delta.
\]
Letting first $\delta\to0$ and then $\epsilon\to0$,
and taking the infimum over $u$, proves the claim since
$
 C_0^\infty(\Omega)\subset W^{1,q}_0(\Omega)
 \subset W^{1,p}_0(\Omega).
$
\end{proof}

%%%%%%%%%%%%%%%%%%%%%%%%%%%%%%%%%%%%%%%%%%%%%%%%%%%%%%%%%%%%%%%%%%%%%%%%%%%%%%%%%%%%%%%%%%
%%%%%%%%%%%%%%%%%%%%%%%%%%%%%%%%%%%%%%%%%%%%%%%%%%%%%%%%%%%%%%%%%%%%%%%%%%%%%%%%%%%%%%%%%%
%%%%%%%%%%%%%%%%%%%%%%%%%%%%%%%%%%%%%%%%%%%%%%%%%%%%%%%%%%%%%%%%%%%%%%%%%%%%%%%%%%%%%%%%%%

\section{The generalized Orlicz case}\label{sect:orlicz}

We now apply the approximation lemma in generalized Orlicz spaces. 
This allows the growth and its dependence on $x$ to be described by a local balance condition, yielding a more general formulation also in the sub-dimensional case.

For a measurable set $E\subset\Omega$, we write
\[
 \phi_E^-(t):=\essinf_{x\in E}\phi(x,t)
 \quad\text{and}\quad
 \phi_E^+(t):=\esssup_{x\in E}\phi(x,t);
\]
we omit the subscript when $E=\Omega$.
We use the \emph{modular} and the \emph{Luxemburg norm}
\[
 \rho_\phi(f):=\int_\Omega\phi(x,|f(x)|)\,dx
 \quad\text{and}\quad
 \|f\|_\phi:=\inf\{\lambda>0:\rho_\phi(f/\lambda)\le1\}.
\]
The space $L^\phi(\Omega)$ consists of the measurable functions
$f$ with $\rho_\phi(f/\lambda)<\infty$ for some $\lambda>0$,
or equivalently, with $\|f\|_\phi<\infty$. 
The Sobolev space $W^{1,\phi}(\Omega)$ consists of the functions in
$W^{1,1}_{\loc}(\Omega)$ for which both $u$ and $|\nabla u|$ belong to
$L^\phi(\Omega)$. We set
$\|u\|_{W^{1,\phi}}:=\|u\|_\phi+\|\nabla u\|_\phi$ and
$\rho_{W^{1,\phi}}(u):=\rho_\phi(u)+\rho_\phi(|\nabla u|)$.
For weak $\Phi$-functions the Luxemburg norms are technically only quasi-norms, 
but since they are equivalent to norms, constants in the triangle inequality are harmless.
We next define the non-weighted condition \azero{} and a local balance condition \aonen{s}. 

\begin{definition}\label{azero}\label{aonen}
We say that $\phi\in \Phiw(\Omega)$ satisfies \azero{} if there exists $\beta_0\in (0, 1]$ 
such that 
\[
\phi^+(\beta_0)\le 1 \le \phi^-(\tfrac1{\beta_0}).
\] 
We say that it satisfies \aonen{s} if
there exists $\beta_1\in(0,1]$ such that
\[
\phi^+_B(\beta_1 t)\le\phi^-_B(t)\qquad\text{whenever }\ t^s\in\big[1,\tfrac1{|B|}\big]
\]
holds in every ball $B\subset\Omega$. If the inequality holds whenever
$t^s\in[0,\frac1{|B|}]$, then we call it \emph{strong \aonen{s}}.
\end{definition}

The idea behind the parameter in \aonen{s} is that we 
work with the norm $\|u\|_\phi$ under the 
a priori assumption $u\in L^s$, cf.\ \cite{HasO23}. 
Adding a linear term extends the comparison condition to small arguments, 
transforming \aonen{s} to strong \aonen{s}.
On bounded sets this does not change the associated space under \azero{}.

\begin{lemma}\label{lem:strongerA1}
Assume that $\phi\in\Phiw(\Omega)$ satisfies \azero{} and set 
$\hat\phi(x,t):=\phi(x,t)+t$. 
If $\phi$ satisfies \aonen{s}, then $\hat\phi$ satisfies strong \aonen{s}. 
\end{lemma}
\begin{proof}
The inequality $\hat\phi^+_B(\beta_1 t)\le \hat\phi^-_B(t)$ follows directly from 
\aonen{s} of $\phi$ when $t^s\in [1, \frac 1{|B|}]$. Assume next 
that $t^s<1$. Choose $\hat\beta_1 := \min\{\beta_1, \frac{\beta_0}{2L}\}$. 
By \ainc{1} and \azero{}, 
\[
\phi^+(\hat\beta_1 t) 
\le
\tfrac {L \hat\beta_1 t}{\beta_0} \phi^+(\beta_0) 
\le \tfrac t2.
\]
Thus $\hat\phi^+(\hat\beta_1 t)\le \frac t2 + \hat\beta_1 t \le t\le \hat\phi^-(t)$, so the required inequality holds 
also in this case.
\end{proof}

Before mollifying, we reduce to bounded functions with compact support.
This reduction uses only cut-off functions and truncation of the values.

\begin{lemma}[Lemma~6.4.2, \cite{HarH19}]\label{lem:dens-compact-orlicz}
Let $\phi\in\Phiw(\Rn)$ satisfy \adec{} and $L^\phi\subset L^1_{\text{loc}}$.
Then, for every $u\in W^{1,\phi}(\Rn)$, there exists a sequence of bounded, compactly supported functions $u_m\in W^{1,\phi}(\Rn)$ such that $\supp u_m\subset\supp u$ and $u_m\to u$ in $W^{1,\phi}(\Rn)$.
\end{lemma}

The next lemma gives a uniform estimate for mollification. Its second
part also controls the error after subtracting a bounded vector field.
This is the estimate needed for strong convergence. Here and below, $\eta_r$ denotes the
standard mollifier fixed in Section~\ref{sect:superdimensional}.

In the following lemmas $C^{0,0}$ and $[u]_0$ are interpreted as $L^\infty$ and $\|u\|_\infty$, 
respectively, corresponding to the case $s=n$ and $\gamma=0$. 

\begin{lemma}\label{lem:bdd-convolution}
Let $\phi\in \Phiw(\Rn)$ satisfy \adec{} and strong \aonen{s}, $s\ge n$, and  
$u\in W^{1,\phi}(\Rn) \cap C^{0,\gamma}(\Rn)$, $\gamma:=1-\frac ns$. 
\begin{enumerate}
\item
$\|\nabla(u*\eta_r)\|_\phi\lesssim \|\nabla u\|_\phi + [u]_\gamma$.
\item
If $r, M>0$ and $[u]_\gamma\le M$, then
\[
\rho_\phi\big((\nabla u - S)*\eta_r\big)\lesssim  \rho_\phi(\nabla u-S)
\]
for every $S\in L^\infty(\Rn; \Rn)$ with $\|S\|_\infty\le M r^{-\frac ns}$. 
The implicit constant depends only on 
$n$, $M$, and the structural constants of $\phi$.
\end{enumerate}
\end{lemma}
\begin{proof}
Since $\phi^-_B$ is comparable to a convex function (cf.\ 
\cite[Lemma 4.3.1]{HarH19}), it follows for $f\in L^1_{\loc}(\R^n)$ that
\begin{equation}\label{eq:jensen}
 \phi^-_B\left(\int_\Rn|f|\eta_r\,dx\right)
 \lesssim
 \int_\Rn\phi^-_B(|f|)\eta_r\,dx.
\end{equation}

(1) Since $\nabla(u*\eta_r) = \nabla \eta_r*(u-u_B)$ for $B:=B_r(x)$, we find that
\[
|\nabla(u*\eta_r)(x)| 
\le 
\frac{c_1(n)}r\fint_{B} |u-u_{B}|\, dy  
\le c_2(n) r^{\gamma-1} [u]_\gamma = c_2(n) r^{-\frac ns} [u]_\gamma
\] 
for some $c_2(n)\ge c_1(n)\ge 1$.

If $\|\nabla u\|_\phi+[u]_\gamma=0$, the assertion is immediate. Otherwise, 
define
$\hat u:=\frac{u}{2\|\nabla u\|_\phi+[u]_\gamma}$.
By the definition of the Luxemburg norm,
we have $\rho_\phi(|\nabla\hat u|)\le 1$. From $[\hat u]_\gamma\le1$ it follows that 
$|\nabla(\hat u * \eta_r)(x)|^s\le \frac{c}{|B|}$ for some $c\ge 1$ depending only on $n$ 
by the estimate above. Thus we can apply strong \aonen{s} at
$t=c^{-1/s}\,|\nabla(\hat u*\eta_r)(x)|$.
We conclude by this, the Jensen inequality \eqref{eq:jensen}, and Fubini's theorem that
\begin{align*}
\int_\Rn \phi(x, \beta_1 c^{-1/s}\,|\nabla(\hat u*\eta_r) (x)|)\, dx
&\le
\int_\Rn \phi_{B_r(x)}^-(|(\nabla\hat u*\eta_r)(x) |) \, dx\\
&\lesssim
\int_\Rn \int_{\Rn} \phi_{B_r(x)}^-(|\nabla\hat u (y) |)\eta_r(x-y)\,dy\, dx \\
&\le
\int_\Rn \phi(y,|\nabla\hat u (y) |) \left[\int_{\Rn} \eta_r(x-y)\,dx\right]\, dy  \\
&= \int_\Rn \phi(y,|\nabla\hat u|)\, dy \le 1.
\end{align*}
In view of the definition of the norm, it follows that 
$\|\nabla \hat u *\eta_r\|_\phi\lesssim 1$. This proves (1).

(2) Let $r, M>0$ with $[u]_\gamma \le M$, and suppose that
$S\in L^\infty(\Rn;\Rn)$ satisfies $\|S\|_\infty\le M r^{-\frac ns}$. 
In view of the earlier estimate, 
\[\begin{split}
|\nabla(u*\eta_r)-S*\eta_r|
& \le |\nabla(u*\eta_r)| + |S*\eta_r| \\
& \le c_2(n) r^{-\frac ns} [u]_\gamma +  \|S\|_\infty \le (M+1)(c_2(n)+1)r^{-\frac ns}
\end{split}\]
Choose $\kappa:= \min\left\{1,\frac{1}{(M+1)(c_2(n)+1)|B_1|^{1/s}}\right\}$.
Then
\[
\bigl(\kappa\,|\nabla(u*\eta_r)-S*\eta_r|\bigr)^s\le \frac1{|B_r|}.
\]
Again, it follows by strong \aonen{s} at $t=\kappa\,|\nabla(u*\eta_r)(x)-S*\eta_r(x)|$, 
the Jensen inequality \eqref{eq:jensen}, and Fubini's theorem that
\begin{align*}
\int_\Rn \phi\big(x, \beta_1\kappa\,\big|\nabla(u*\eta_r)(x)-S*\eta_r(x)\big|\big)\, dx
&\le
 \int_\Rn \phi_{B_r(x)}^-\big(\big|(\nabla u-S)*\eta_r(x)\big|\big)\, dx\\
&\lesssim
\int_\Rn\int_\Rn \phi_{B_r(x)}^-(|\nabla u(y)-S(y)|)\eta_r(x-y)\,dy\,dx\\
&\le
\int_\Rn \phi(x,|\nabla u(x)-S(x)|)\, dx. 
\end{align*}
The conclusion follows from \adec{}, which allows us to remove the factor $\beta_1\kappa$ on the left-hand side.
\end{proof}

We next prove the analogue of Lemma~\ref{lem:mollification} with the local balance condition.

\begin{lemma}\label{lem:orlicz-holder}
Let $\phi\in \Phiw(\Rn)$ satisfy \adec{} and \aonen{s}, $s\ge n$.
If $u\in W^{1,\phi}(\Rn)\cap C^{0,\gamma}(\Rn)$, $\gamma:=1-\frac ns$, has compact support, then
$u*\eta_r\to u$ in $W^{1,\phi}(\Rn)$ as $r\to0$.
\end{lemma}
\begin{proof}
Let $0<\epsilon<1$ and define
\[
 \Omega:=\{x\in\Rn:\dist(x,\supp u)<1\}
 \quad\text{and}\quad
 U:=\{x\in\Rn:\dist(x,\supp u)<2\}.
\]
Since $\Omega$ is bounded, it follows from Lemma~\ref{lem:strongerA1} and 
\cite[Theorem~3.2.6]{HarH19} that
$\|f\|_\phi\approx\|f\|_{\hat\phi}$ for every $f\in L^\phi(\R^n)$ supported in
$\Omega$, where $\hat\phi(x,t):=\phi(x,t)+t$.
By Lemma~\ref{lem:strongerA1}, applied on $U$, $\hat\phi$ satisfies
strong \aonen{s} on $U$; it also satisfies \adec{q} with the same exponent as $\phi$.

Applying \cite[Theorem~3.5.2]{HarH19} to each component of $\nabla u$, 
we obtain $S\in L^\infty(\Rn;\Rn)$ with 
$\supp S \subset  \supp u \subset \Omega$ such that $\|\nabla u-S\|_{\hat\phi}\le\epsilon$
(we can for instance use $S:=\nabla u \chi_{\{|\nabla u|\le m\}}$ for large $m$). 
Define $u_r:=u*\eta_r$ and $M:=[u]_\gamma$ and 
consider $r\in (0,1)$ such that $\|S\|_\infty \le M r^{-\frac ns}$.
For centers $x\in\Omega$, all these balls are contained in $U$.
By Lemma~\ref{lem:bdd-convolution}(2) and the modular--norm inequality \cite[Lemma 3.2.9]{HarH19},
\begin{align*}
 \|(\nabla u-S)*\eta_r\|_\phi
 &\approx 
 \|(\nabla u-S)*\eta_r\|_{\hat\phi} \\
 &\lesssim
 \max\Big\{\rho_{\hat\phi}((\nabla u-S)*\eta_r),
 \rho_{\hat\phi}((\nabla u-S)*\eta_r)^{\frac1q}\Big\} \\
  &\lesssim
 \max\Big\{\rho_{\hat\phi}(\nabla u-S),
 \rho_{\hat\phi}(\nabla u-S)^{\frac1q}\Big\} \\
 &\lesssim\|\nabla u-S\|_{\hat\phi}^{1/q}
 \le \epsilon^{1/q},
\end{align*}
where the implicit constants may depend on $M$ and the local structure constants on $U$,
but are independent of $\epsilon$, $S$, and $r$.
Thus
\begin{align*}
 \|\nabla u_r-\nabla u\|_\phi
 &\lesssim
 \|(\nabla u-S)*\eta_r\|_\phi
 +\|S*\eta_r-S\|_\phi
 +\|\nabla u-S\|_\phi\\
 &\le C\epsilon^{1/q}+\|S*\eta_r-S\|_\phi.
\end{align*}
Since $\supp S\subset\Omega$ and $S\in L^q(\Omega)$, it follows that $S*\eta_r\to S$ in 
$L^q(\Omega)$. Since $L^q(\Omega)$ embeds continuously into $L^\phi(\Omega)$ by 
\cite[Lemma~3.7.7]{HarH19}, $\|S*\eta_r-S\|_\phi\to 0$.
Thus, first letting $r\to 0$ and then $\epsilon\to 0$ in the
previous estimate, we obtain that $\nabla u_r\to\nabla u$ in $L^\phi(\Rn)$.
Since also $u\in L^\infty(\Rn)$, we have
$u_r\to u$ in $L^q(\Omega)\subset L^\phi(\Omega)$
and hence in $L^\phi(\Rn)$, which completes the proof.
\end{proof}

We conclude the section with an improved density result for 
generalized Orlicz spaces. We use the boundedness obtained in Lemma~\ref{lem:dens-compact-orlicz}
in the sub-dimensional case while the 
improved H\"older exponent of Lemma~\ref{lem:key} is used in the super-dimensional case.

\begin{theorem}[Density for generalized Orlicz--Sobolev space]\label{thm:orlicz-super}
Let $1\le p\le q<\infty$ and $s:=(\max\{p, n\}-1)n/(n-1)$.
Suppose that $\phi\in\Phiw(\Rn)$ satisfies \ainc{p},
\adec{q}, and \aonen{s}.
Then $C_0^\infty(\Rn)$ is dense in $W^{1,\phi}(\Rn)$.
\end{theorem}
\begin{proof}
First, \aonen{s} and \adec{q} imply that $\phi$ satisfies
\azero{} on every bounded set, which implies $L^\phi(\R^n)\subset L^1_{\loc}(\R^n)$ \cite[Lemma~3.7.7]{HarH19}. 
By Lemma~\ref{lem:dens-compact-orlicz} and \adec{q}, bounded functions
with compact support are dense in $W^{1,\phi}(\Rn)$. Thus it suffices
to show that every such function can be approximated by smooth functions
in $W^{1,\phi}(\Rn)$.

Choose a bounded open set $U$ containing $\supp u$.
Then \azero{} and \ainc{p} imply $t^p\lesssim\phi(x,t)+1$ on $U$, and hence
$u\in W^{1,p}(\Rn)$ \cite[Lemma~3.7.7]{HarH19}. 
Moreover, \adec{q} and $\|u\|_{W^{1,\phi}}<\infty$ imply $\rho_\phi(|\nabla u|)<\infty$ 
\cite[Lemma~3.2.9]{HarH19}.
Define $\gamma_*:=1-\frac ns=1-\frac{n-1}{p-1}$ when $p>n$. 
From Lemma~\ref{lem:key} and dominated convergence we obtain compactly
supported functions $v_j\in W^{1,p}(\Rn)\cap C^{0,\gamma_*}(\Rn)$ with
$|\nabla v_j|\le|\nabla u|$ and $\supp v_j \subset  \supp u$, and
\[
 \|u-v_j\|_\infty+\rho_\phi(|\nabla u-\nabla v_j|)\longrightarrow 0.
\]
It thus suffices to approximate functions in 
$W^{1,\phi}(U)\cap C^{0,\gamma_*}(U)$ by smooth functions.
When $p\le n$, this step is not needed and we instead approximate 
$W^{1,\phi}(U)\cap L^\infty(U)$-functions. 
In either case, the convolutions give the desired approximations, by Lemma~\ref{lem:orlicz-holder}.
\end{proof}

%%%%%%%%%%%%%%%%%%%%%%%%%%%%%%%%%%%%%%%%%%%%%%%%%%%%%%%%%%%%%%%%%%%%%%%%%%%%%%%%%%%%%%%%%%
%%%%%%%%%%%%%%%%%%%%%%%%%%%%%%%%%%%%%%%%%%%%%%%%%%%%%%%%%%%%%%%%%%%%%%%%%%%%%%%%%%%%%%%%%%
%%%%%%%%%%%%%%%%%%%%%%%%%%%%%%%%%%%%%%%%%%%%%%%%%%%%%%%%%%%%%%%%%%%%%%%%%%%%%%%%%%%%%%%%%%
\section{Improved conditions for double phase maximal regularity}
\label{sect:maximal}

We apply the density result to scalar weak subsolutions of
quasilinear equations with double phase growth.
Density allows us to use truncated functions as tests.
The resulting flattening estimate, combined with Morrey's
inequality, yields improved H\"older continuity of subsolutions in the 
super-dimensional case.

\begin{proposition}\label{thm:radial-subsolution}
Let $n<p\le q$, $0<\alpha\le1$, and
$0\le a\in C^{0,\alpha}_{\loc}(\Omega).$ Suppose that
\[
H(x,t):=t^p+a(x)t^q
\quad\text{and}\quad
q-p\le\alpha\frac{p-1}{n-1}.
\]
Let $\Lambda,\nu>0$ and suppose $A:\Omega\times\Rn\to\Rn$ with $A(x,0)=0$ satisfies
\[
\tfrac1\Lambda |A(x,\xi)| \,|\xi| 
\le 
H(x,|\xi|)
\le
\tfrac1\nu A(x,\xi)\cdot\xi 
\]
for every $x\in \Omega$ and $\xi\in\Rn$.
If $u\in W^{1,p}_{\loc}(\Omega)$ satisfies
$H(\cdot,|\nabla u|)\in L^1_{\loc}(\Omega)$ and
\begin{equation}\label{eq:radial-weak-subsolution}
 \int_\Omega A(x,\nabla u)\cdot\nabla\zeta\,dx\le0
\qquad \text{for every }\ \zeta\in C_0^\infty(\Omega) \ \ \text{with}\ \ \zeta\ge 0,
\end{equation}
then $u\in C^{0,\gamma_*}_{\loc}(\Omega)$, where
$\gamma_*= \frac{p-n}{p-1}$.
\end{proposition}

\begin{proof}
By Morrey's embedding, $u\in C^{0,\gamma}_{\loc}(\Omega)$ with $\gamma:=1-\frac np<\gamma_*$. 
Let $B_{2R}\Subset \Omega$ and suppose the cut-off function $\zeta\in C_0^\infty(B_{2R})$ 
satisfies $0\le\zeta\le 1$ with $\zeta=1$ on $B_R$.
For $s<t$ set
\[
 T_{s,t}(\tau):=\min\{\max\{\tau-s,0\},t-s\}
 \quad\text{and}\quad
 \tilde\zeta:=\zeta T_{s,t}(u).
\]
Since $\tilde\zeta\ge 0$, $\supp \tilde\zeta\Subset B_{2R}$, and $\mathcal F_H(\tilde\zeta,B_{2R})<\infty$, by the construction of smooth functions in the proof of Theorem~\ref{thm:pq}, one can find
$\zeta_j\in C_0^\infty(\Omega)$, $j\in \mathbb N$, with $\zeta_j\ge0$ and $\supp \zeta_j\Subset B_{2R}$ such that
\[
 \int_{B_{2R}} H(x,|\nabla\zeta_j-\nabla\tilde\zeta|)\,dx \ \longrightarrow\ 0
 \quad \text{as }\ j\to \infty.
\]
This together with H\"older's inequality implies
\[
\begin{split}
 \int_{B_{2R}}|A(x,\nabla u)|\,|\nabla\zeta_j-\nabla\tilde\zeta|\,dx
 &\le
 \Lambda\left(\int_{B_{2R}}|\nabla u|^p\,dx\right)^{\frac{p-1}{p}}
  \left(\int_{B_{2R}}|\nabla\zeta_j-\nabla\tilde\zeta|^p\,dx\right)^{\frac1p}\\
 &\quad+
 \Lambda\left(\int_{B_{2R}}a(x)|\nabla u|^q\,dx\right)^{\frac{q-1}{q}}
  \left(\int_{B_{2R}} a(x)|\nabla\zeta_j-\nabla\tilde\zeta|^q\,dx\right)^{\frac1q}\\
 &\ \longrightarrow\ 0
\quad \text{as }\ j\to \infty.
\end{split}
\]
Since \eqref{eq:radial-weak-subsolution} holds for every $\zeta_j$, we obtain it for $\tilde\zeta:=\zeta T_{s,t}(u)$
at the limit, i.e. 
\[
\int_\Omega A(x,\nabla u)\cdot\nabla(\zeta T_{s, t}(u))\,dx\le 0.
\]

For any ball $B_{r}\subset B_{R}$, let $s:=\inf_{B_r}u$ and $t:=\sup_{B_r}u$. Suppose $t-s>0$. 
We estimate $\nu\, |\xi|^p \le A(x,\xi)\cdot\xi$, note that $\chi_{\{s \le u \le t\}} \zeta=1$ in 
$B_r$ and that $\nabla (T_{s,t}(u))=\chi_{\{s\le u \le t\}} \nabla u$ by the chain rule, 
and use the conclusion of the previous paragraph to obtain
\[\begin{split}
\nu \int_{B_r}
|\nabla u|^p\,dx
&\le
\int_{B_{2R}}
\chi_{\{s\le u\le t\}} \zeta A(x, \nabla u)\cdot \nabla u\,dx
\le
 -\int_{B_{2R}} T_{s,t}(u)A(x,\nabla u)\cdot\nabla\zeta\,dx\\
& \le (t-s) \Lambda \|\nabla \zeta\|_{\infty} \int_{B_{2R}} |\nabla u|^{p-1}+a(x)|\nabla u|^{q-1}\,dx 
=: (t-s) M.
\end{split}\]
Therefore, by the Morrey inequality \eqref{eq:morrey-oscillation},
\[
 (t-s)^p
 \le c r^{p-n}\int_{B_r}|\nabla u|^p\,dx
 \le c M r^{p-n}(t-s) \quad \Longrightarrow \quad  \osc_{B_r}u = t-s \le c M^{\frac{1}{p-1}}r^{\frac{p-n}{p-1}}.
\]
This proves the assertion since $B_r\subset B_R$ is arbitrary.
\end{proof}

For local minimizers of the functional $\mathcal F_H$, the improved H\"older continuity allows us
to apply the regularity theorem of Baroni, Colombo and
Mingione~\cite[Theorem~1.3]{BarCM18}.

\begin{corollary}[Gradient regularity for double phase minimizers]
\label{cor:BCM-super}
Let $n<p<q$, $0<\alpha\le1$, and
$0\le a\in C^{0,\alpha}_{\loc}(\Omega)$. Suppose that
\begin{equation}\label{eq:BCM-improved-range}
 q-p<\alpha\frac{p-1}{n-1}.
\end{equation}
If $u\in W^{1,p}_{\loc}(\Omega)$ with $H(\cdot,|\nabla u|)\in L^1_{\loc}(\Omega)$ is a local minimizer of
\[
 v\mapsto\int_\Omega H(x,|\nabla v|)\,dx,
\]
then
$u\in C^{1,\beta}_{\loc}(\Omega)$ for some $\beta\in(0,1)$.
\end{corollary}

\begin{proof}
The minimality inequality shows that
$u$ satisfies the Euler--Lagrange equation
\[
 \int_\Omega
 \big(p|\nabla u|^{p-2}\nabla u
       +qa(x)|\nabla u|^{q-2}\nabla u\big)
 \cdot\nabla\zeta\,dx=0
 \quad \text{for every }\ \zeta\in C_0^\infty(\Omega).
\]
Proposition~\ref{thm:radial-subsolution} therefore yields $u\in C^{0,\gamma_*}_{\loc}(\Omega)$.
Since \eqref{eq:BCM-improved-range} is equivalent to
\[
 q<p+\frac{\alpha}{1-\gamma_*},
\]
the conclusion follows from
\cite[Theorem~1.3]{BarCM18}.
\end{proof}

\section*{Acknowledgments}
J. Ok was supported by the National Research Foundation of Korea (NRF) grants funded by the Korean government (MSIT) (NRF-2022R1C1C1004523 and RS-2025-24533680).

\medskip
\noindent\textbf{Data availability.}
No datasets were generated or analyzed during this study.

\medskip
\noindent\textbf{Conflict of interest.}
The authors declare no conflicts of interest.

\medskip
\noindent\textbf{AI statement.}
ChatGPT was used to checking intermediate calculations, refining parts of the proofs, 
and polishing the language. The final text is fully revised and checked by the authors.

%%%%%%%%%%%%%%%%%%%%%%%%%%%%%%%%%%%%%%%%%%%%%%%%%%%%%%%%%%%%%%%%%%%%%%%%%%%%%%%%%%%%%%%%%%
%%%%%%%%%%%%%%%%%%%%%%%%%%%%%%%%%%%%%%%%%%%%%%%%%%%%%%%%%%%%%%%%%%%%%%%%%%%%%%%%%%%%%%%%%%
%%%%%%%%%%%%%%%%%%%%%%%%%%%%%%%%%%%%%%%%%%%%%%%%%%%%%%%%%%%%%%%%%%%%%%%%%%%%%%%%%%%%%%%%%%

\end{document}